\documentclass[12pt,oneside]{amsart}

      \usepackage{amssymb}
      \usepackage[english]{babel}
      \usepackage{graphicx}
      \usepackage[utf8]{inputenc}
      \usepackage{float}
      
      \usepackage{wrapfig}

      \usepackage[a4paper=true,%
      breaklinks=true,%
      colorlinks=true,%
      linkcolor={blue},%
      citecolor={red},%
      pdftitle={Complexity of natural numbers},%
      pdfauthor={J. Arias de Reyna},%
      pdfkeywords={snapshot},%
      pdfstartview=FitH,
]{hyperref}

      \theoremstyle{plain}
      \newtheorem{theorem}{Theorem}
      \newtheorem{proposition}[theorem]{Proposition}
      
      \newtheorem{corollary}[theorem]{Corollary}
      \newtheorem*{problem}{Problem}
      \newtheorem{conjecture}{Conjecture}
      
      \theoremstyle{definition}
      \newtheorem{definition}[theorem]{Definition}

      \theoremstyle{remark}
      \newtheorem{remark}{Remark}

      \newcommand{\N}{{\mathbb N}}
      \newcommand{\R}{{\mathbb R}}
      
      \newcommand{\card}{\mathop{\rm card }\limits}
      
      \newfont{\cmbsy}{cmbsy10}
      \newfont{\cmmib}{cmmib10}
      \newcommand{\Prob}{\mathop{\Bbb P}}

      \def\PNP{\text{\bf P}\buildrel{?}\over{=}\text{\bf NP}}

      \makeatletter
      \def\@setcopyright{}
      \def\serieslogo@{}
      \makeatother
   
\begin{document}

%


   
   \author{Juan Arias de Reyna$\,\,$}

   \address{Facultad de Matem\'aticas, 
   Universidad de Sevilla, \newline
   Apdo.~1160, 41080-Sevilla, Spain}
   \email{arias@us.es}

%
%



   \title
   {Complexity of natural numbers}






   \dedicatory{English version of the paper: \emph{Complejidad de los números naturales}, Gaceta de la Real Sociedad Matemática Española \textbf{3} (2000) 230--250.}

   \date{May, 2000}


   \maketitle



   \section{Introduction.}
   
   \subsection{Complexity of a natural number.}
   Our purpose is to explore what seems to be a trivial question that may be 
   understood by  highschool students in their early teens,
   but with very deep relationships. 
   
   Lately 
   I have been interested in one of the mathematical problems that I consider most 
   important: the $\PNP$ problem. In this case the first difficulty is 
   to explain the problem to a professional mathematician, say to an expert
   in Analysis. This is not a minor issue, I think that the problem $\PNP$
   may be put as an inequality. Hence to explain the question adequately, 
   so that it is understood by an expert in Analysis, maybe
   the first step in the solution of the problem. 
   
   The question I 
   shall discuss here arose
   while trying to obtain this explanation.
   
   We start with the main question: Given a natural number $n$, how many $1$'s
   are needed to write $n$? For example
   \begin{displaymath}
   19=1+(1+1)(1+1+1)(1+1+1)
   \end{displaymath}
   so that nine $1$'s suffice to write $19$. We shall say that the complexity
   of $19$ is less than or equal to $9$, and we shall write this as 
   $\Vert 19\Vert\le 9$. Of course, the complexity of $19$ will be  the 
   number of $1$'s  in the most economical representation of $19$. We only 
   admit expressions with sums and products.
   
   The first values of the complexity function may be easily computed
   \begin{displaymath}
   1,2,3,4,5,5,6,6,6,7,8,7,8,8,8,8,9,8,9,9,\dots
   \end{displaymath}
   We see that this is not a monotonic sequence: $8=\Vert11\Vert>\Vert12\Vert=7$.
   
   When in our investigations we find any sequence of natural numbers, there is 
   something we must do: look in \href{http://oeis.org}
   {\emph{The On-Line Encyclopedia of Integer 
   Sequences}} of Sloane and Plouffe \cite{SP}. In it we find this sequence and
   a reference to a paper by Guy \cite{G} where it is defined and analyzed.
   
   \section{Complexity of a natural number.}
   We have defined the complexity as a function $n\mapsto\Vert n\Vert$ of $\N\to\N$
   such that for every pair of natural numbers $m$ and $n$ we have
   \begin{displaymath}
   \Vert 1\Vert=1, \qquad\Vert m+n\Vert\le \Vert m\Vert +\Vert n\Vert,\qquad
   \Vert m\cdot n\Vert\le \Vert m\Vert +\Vert n\Vert.
   \end{displaymath}
   In fact it is the largest function satisfying these conditions. To prove this
   and other assertions it is useful to introduce  the concept of \emph{expression}.
   
   \subsection{Definition of expression.}
   An expression is a sequence of symbols. The allowed symbols are \texttt{
   x, +, (, )}. Not every sequence of these symbols is an expression. Examples
   of expressions are:
   \begin{displaymath}
   \texttt{(x + x)}; \quad \texttt{(x+(xx))};\qquad
   \texttt{(x+((x+x)((x+(x+x))(x+(x+x)))))}.
   \end{displaymath}
   The formal definition is inductive:
   \begin{itemize}
   \item[(a)] \texttt{x} is an expression.
   \item[(b)] If \texttt{A} and \texttt{B} are expressions, then 
   \texttt{(A+B)} and \texttt{(AB)} are also expressions.
   \item[(c)] The only expressions are those obtained by  repeated 
   applications of rules (a) and (b). 
   \end{itemize}
   We define the value of an expression \texttt{A} as the number $v(\texttt{A})$ 
   that results when replacing \texttt{x} by $1$.
   Again we use induction to define the value of an expression: $v(\texttt{x})=1$, 
   and if \texttt{A} and \texttt{B} are expressions then $v(\texttt{(A+B)})=v(
   \texttt{A})+v(\texttt{B})$ and $v(\texttt{(AB)})=v(
   \texttt{A})v(\texttt{B})$.
   
   Given an expression we may define its complexity as the number of letters 
   \texttt{x}
   it contains, for example $\Vert\texttt{(x+(xx))}\Vert=3$. Let $\mathcal{E}$
   be the set of expressions. We may translate the definition of the complexity 
   as 
   \begin{displaymath}
   \Vert n\Vert=\inf\{\Vert\texttt{A}\Vert: \texttt{A}\in\mathcal{E} \;
   \textrm{and}\; v(\texttt{A})=n\}.
   \end{displaymath}
   
   If we want to compute the value of $\Vert n\Vert$ we may use the following
   Proposition.
   \begin{proposition}\label{computing}
   For each natural number $n>1$
   \begin{displaymath}
   \Vert n\Vert=\min_{\substack{2\le d\le \sqrt{n},\,\, d\mid n\\
   1\le j\le n/2}}\left\{\Vert d\Vert+\Vert n/d\Vert,\quad \Vert j\Vert+\Vert
   n-j\Vert\right\}
   \end{displaymath}
   \end{proposition}
   
   \begin{proof}
   Let \texttt{E} an optimal expression for $n$, i.~e.~one that gives its complexity
   $\Vert n\Vert =\Vert \texttt{E}\Vert$. As an expression that is not 
   \texttt{x} we will
   have \texttt{E}=\texttt{(A+B)} or \texttt{E}=\texttt{(AB)}. Let $a=v(\texttt{A})$
    and $b=v( \texttt{B})$. Then either $n=a+b$ and $\Vert n\Vert=\Vert a\Vert+
    \Vert b\Vert$ or $n=ab$ and $\Vert n\Vert=\Vert a\Vert+
    \Vert b\Vert$.  In the first case if $j$ is the least of $a$ and $b$  we will
    have $1\le j\le n/2$, and in the second case  if $d$ is the least of $a$ and $b$,
    then $d$ will be a divisor of $n$ with $2\le d\le \sqrt{n}$.  Of course for the 
    reasoning to be valid  we must check that if \texttt{E} is an optimal expression
    for $n$, then \texttt{A} and \texttt{B} must be optimal expressions for $a$ and 
    $b$ respectively. We leave this check to the reader.
    \end{proof}
    
    Using the above Proposition and the mathematical software Mathematica we have
    computed the values of $\Vert n\Vert$ for $1\le n\le 200\,000$.    
    
    \section{Bounds.}
    \begin{proposition}
    Let $P\colon\N\to\R$    be a function satisfying
    \begin{displaymath}
    P(1)=1,\quad P(n+m)\le P(n)+P(m),\quad P(n\cdot m)\le P(n)+P(m).
    \end{displaymath}
    Then for each $n\in\N$ we have $P(n)\le\Vert n\Vert$.
    \end{proposition}
    \begin{proof}
    It is easy to see by induction that for each expression \texttt{A}, 
    we have $P(v(\texttt{A}))
    \le \Vert \texttt{A}\Vert$.
    It is true for $\texttt{A}=\texttt{x}$, and, if it is true for \texttt{A}
    and \texttt{B}  then it is true for \texttt{(A+B)} and \texttt{(AB)}. For 
    example, for the product:
    \begin{displaymath}
    P\bigl(v(\texttt{(AB)})\bigr)=P\bigl(v(\texttt{A})v(\texttt{B})\bigr)
    \le P\bigl(v(\texttt{A})\bigr)+P\bigl(v(\texttt{B})\bigr)\le
    \Vert \texttt{A}\Vert + \Vert \texttt{B}\Vert =\Vert \texttt{(AB)}\Vert,
    \end{displaymath}
    and a similar argument is valid for the sum. (Observe that by the definition
    of $v$ we have $v(\texttt{(A+B)})=v(\texttt{A})+v(\texttt{B})$ and 
    $v(\texttt{(AB)})=v(\texttt{A})v(\texttt{B})$).
    
    Now in $P(v(\texttt{A}))\le \Vert A\Vert$ we take the minimum  over all 
    expressions
    \texttt{A} such that $n=v(\texttt{A})$. In this way we get $P(n)\le \Vert n\Vert$.
    \end{proof}
    
    \begin{corollary}
    For each natural number $n$ we have $\log_2(1+n)\le \Vert n\Vert$.
    \end{corollary}
    \begin{proof}
    It is sufficient to check the properties of $P(n)=\log_2(1+n)$.
    \end{proof}
    Later, in Corollary \ref{antnueve}, we will obtain a better inequality.
    
    \subsection{Upper bounds.}
    Now we get an upper bound. To this end we define a new function 
    $L\colon\N\to\N$.
    
    \begin{definition}
    We define the function $L$ inductively:
    \begin{itemize}
    \item[(a)] $L(1)=1$. 
    \item[(b)] If $p$ is a prime number, then $L(p)=1+L(p-1)$.
    \item[(c)] If $n=p_1p_2\cdots p_k$ is a product of prime numbers (may be
    repeated), then $L(p_1p_2\cdots p_k)=L(p_1)+L(p_2)+\cdots + L(p_k)$.
    \end{itemize}
    \end{definition}
    It is clear from this definition that if $n=ab$ with $a$ and $b\ge2$ then we
    will have $L(ab)=L(a)+L(b)$.
    
    \begin{proposition}
    For each $n\in\N$ we have 
    \begin{displaymath}
    \Vert n\Vert\le L(n).
    \end{displaymath}
    \end{proposition}
    
    \begin{proof}
    We may prove this by induction. For $n=1$ we have $\Vert 1\Vert=L(1)=1$. Assume
    that $\Vert k\Vert\le L(k)$ for each $k<n$.
    There are two possibilities: if $n=p$ is a prime number
    \begin{displaymath}
    \Vert p\Vert \le \Vert p-1\Vert +\Vert 1\Vert =\Vert p-1\Vert +1\le L(p-1)+1=L(p).
    \end{displaymath}
    If  $n$ is composite $n=ab$ with $a$ and $b>2$, 
    \begin{displaymath}
    \Vert n\Vert \le \Vert a\Vert +\Vert b\Vert \le L(a)+L(b)=L(ab)=L(n).
    \end{displaymath}
    \end{proof}
    
    \begin{proposition}
    For each $n\ge2$ we have
    \begin{displaymath}
    L(n)\le\frac{3}{\log 2}(\log n).
    \end{displaymath}
    \end{proposition}
    
    \begin{proof}
    Since $L(2)=2$ and $L(3)=3$ the result is true for  $n = 2$  and  $n = 3$.
    
    Assume now that $n>3$ and that the Proposition is true for all natural 
    numbers strictly less than $n$.
    
    If $n=p$ is a prime number we have
    \begin{equation}\label{eq1}
    L(p)=1+L(p-1)=1+2+L\Bigl(\frac{p-1}{2}\Bigr)\le 3+\frac{3}{\log2}
    \log\Bigl(\frac{p-1}{2}\Bigr).
    \end{equation}
    We want this to be 
    \begin{displaymath}
    \le\frac{3}{\log2}(\log p).
    \end{displaymath}
    Hence we must check that 
    \begin{equation}\label{eq2}
    3\le \frac{3}{\log 2}\log\Bigl(\frac{2p}{p-1}\Bigr),
    \end{equation}
    which is easily proved for $p\ge3$.
    
    If $n=ab$ with $a$ and $b\ge2$, we have
    \begin{displaymath}
    L(ab)=L(a)+L(b)\le\frac{3}{\log 2}(\log a)+\frac{3}{\log 2}(\log b)=
    \frac{3}{\log 2}(\log ab).
    \end{displaymath}
    \end{proof}

    \begin{remark}
    We do not know if the constant $3/\log2$ in the above theorem is optimal.
    The proof makes one suspect that the quotient $L(n)/\log n$ may be large
    when $n=p_k$ is a prime such that there exists a sequence of primes 
    $(p_j)_{j=1}^k$ with $p_{j+1}=2p_{j}+1$. For example,  $89$,
    $179$, $359$, $719$, $1439$, $2879$ is such a sequence of prime numbers, 
    and the maximum value of the quotient $L(n)/\log n$ that we know is
    \begin{displaymath}
    \frac{L(2879)}{\log 2879}=  3.766384578\dots  < 4.328085123\dots   
    =\frac{3}{\log 2}.
    \end{displaymath}
    The main difference between the two functions $L(\cdot)$ and $\Vert\cdot\Vert$
    is that $L(\cdot)$ is additive and $\Vert\cdot\Vert$ is not. For each pair
    of numbers $n$ and $m$ greater than $1$ we have $L(mn)=L(m)+L(n)$. On the other
    hand there exist pairs $n$, $m$ of numbers greater than $1$ and such that
    $\Vert mn\Vert<\Vert m\Vert+\Vert n\Vert$. In such a case we shall say 
    that $n\cdot m$ is a 
    bad factorization.
    
    In figure 1 we put a dot  at each point $(n,m)$ such that $n\cdot m$ is a bad
    factorization. The figure contains all the factors $n$ and $m\le 60$.

    \begin{figure}[H]
    \begin{center}
    \includegraphics[width = \textwidth]{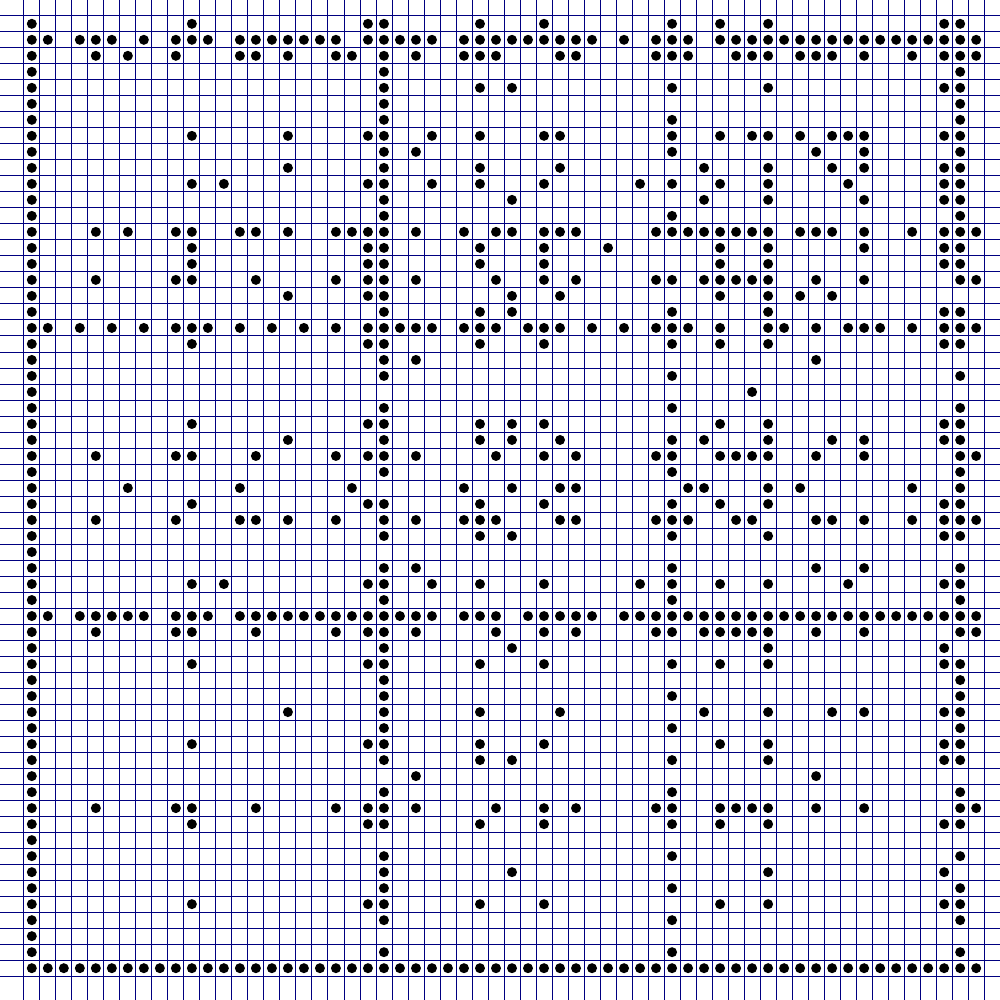}
   \caption{Bad Factors.}
    \end{center}
    \label{tres}
    \end{figure}  
    
    $1\cdot n$ is always a bad factorization. In the figure we see some 
    other surprising
    regularities. There are some conspicuous (vertical and horizontal) aligned points.
    Especially note the verticals at $n=23$, $41$, $59$, which deserve an explanation.
    
    These numbers, we may call them \emph{bad factors},  appear to have  
    great complexity. We define the \emph{number with great complexity}   $n_k$
    as the number $n_k$ that is the less solution to $\Vert n\Vert =k$. The
    first values of this sequence are
    \begin{gather*}
    1,  \quad 2,  \quad 3,  \quad 4,  \quad 5,  \quad 7,  \quad 10,  \quad 11,  
    \quad 17,  \quad 22,  \quad 23,  \quad 41,  \quad 47,  \quad 59,  \\
     89, 
    \quad 107,  \quad 167, \quad 179,  \quad 263,  \quad 347,  \quad 467 ,  \quad 683,  
    \quad 719,  \quad  1223, \\ 1438,  \quad 1439, \quad 2879,  \quad 3767,  
    \quad 4283,\quad 6299,  \quad 10079,  \quad 11807, \\
    15287,  \quad
    21599,
    \quad 33599,  \quad \dots 
    \end{gather*}
    This sequence appears in \cite{SP} with some errata. In this way 
    we find the reference to Rawsthorne \cite{R}.
    
    \end{remark}
    
    \section{Mean values.}
    There is another proof of $\Vert n\Vert\le 3\log n/\log2$. We observe that 
    if we write $n$ in binary $n=\sum_{j=0}^{k-1}\varepsilon_j 2^j+2^{k}$ we have
    a means to express $n$:
    \begin{displaymath}
    n=\varepsilon_0+2(\varepsilon_{1}+2(\varepsilon_{2}+\cdots+
    2(\varepsilon_{k-2}+2(\varepsilon_{k-1}+2))\cdots)).
    \end{displaymath}
    If we substitute each $2$ by $1+1$ and observe that each $\varepsilon_j$
    is equal to $0$ or $1$, we have an expression for $n$ that uses at most  
    $2k+k$ ones, and where $k$ is determined by $2^k\le n<2^{k+1}$. It follows
    that $\Vert n\Vert\le 3\log n/\log2$. 
    
    The above reasoning proves that the function $L_2(n)=2k+\varepsilon_0+
    \varepsilon_1+\cdots+\varepsilon_{k-1}$ is another upper bound for $\Vert n\Vert$.
    The relation between $L_2(n)$ and $L(n)$ is not very simple. 
    Amongst the first 1000 numbers
    we generally have $L(n)\le L_2(n)$ but this inequality has exceptions. 
    The first one is
    $L_2(161)=16<17=L(161)$.  In this range the difference is small.
    
    The function $L_2(n)$ allows us to obtain information about $\Vert \cdot\Vert$. 
    Consider the numbers $n$ that in binary take the form $1\varepsilon_{k-1}\cdots
    \varepsilon_0$, i.~e.~numbers than in  binary have $k+1$ digits.  
    By the above expression
    we have
    \begin{displaymath}
    \Vert n\Vert\le 2k+\varepsilon_0+\cdots +\varepsilon_{k-1}.
    \end{displaymath}
    We may suppose that the $\varepsilon_k$ are independent random variables with 
    mean $1/2$. The inequality of Chernoff (see \cite{C} or \cite{AS} for a simple
    exposition) says that
    \begin{displaymath}
    \Prob\left(\left|\sum\varepsilon_j-k/2\right|<x\sqrt{k}\right)\ge
    1-2e^{-2x^2}.  
    \end{displaymath}
    It follows that $\Prob(\Vert n\Vert \le 2k+k/2+x\sqrt{k})\ge1-2e^{-2x^2}$, 
    and taking $x=\sqrt{\log k}$ we get
    \begin{displaymath}
    \Prob\left(\Vert n\Vert >5k/2+\sqrt{k\log k}\right)\le 2 k^{-2}.
    \end{displaymath}
    Hence between the $2^k$ values of $n$ with $2^k\le n<2^{k+1}$ at most
    $(2/k^2)2^k$ satisfy $\Vert n\Vert >5k/2+\sqrt{k\log k}$. The other ones, most of 
    them, satisfy
    \begin{displaymath}
    \Vert n\Vert \le\frac{5k}{2}+\sqrt{k\log k}=
    \frac{5}{2}\frac{\log n}{\log 2}+ O(\sqrt{\log n\log\log n}).
    \end{displaymath}
    Therefore, for almost all large values of $n$ we have
    \begin{displaymath}
    \Vert n\Vert \le\frac{5}{2}\frac{\log n}{\log 2}+ O(\sqrt{\log n\log\log n}).
    \end{displaymath}
    
    The upper bound $L(n)$ is very good for small values of $n$. For example 
    for the first 220 values  of $n$, $L(n)=\Vert n\Vert$, except for the 
    values in the following table:
    \bigskip
    
\begin{center}    
\begin{tabular}{|c|c|c|}
\hline
$n$ & $\Vert n\Vert$ & $L(n)$\\ \hline\hline
46 & 12 & 13 \\ \hline
47 & 13 & 14 \\ \hline
55 & 12 & 13 \\ \hline
82 & 13 & 14  \\ \hline
83 & 14 & 15  \\ \hline
92 & 14 & 15 \\ \hline
94 & 15 & 16  \\ \hline
110 & 14 & 15 \\ \hline
\end{tabular}
\quad
\begin{tabular}{|c|c|c|}
\hline
$n$ & $\Vert n\Vert$ & $L(n)$\\ \hline\hline
115 & 15 & 16 \\ \hline
118 & 15 & 16  \\ \hline
121 & 15 & 16 \\ \hline
138 & 15 & 16  \\ \hline
139 & 16 & 17  \\ \hline
141 & 16 & 17 \\ \hline
145 & 15 & 16  \\ \hline
161 & 16 & 17 \\ \hline
\end{tabular}
\quad
\begin{tabular}{|c|c|c|}
\hline
$n$ & $\Vert n\Vert$ & $L(n)$\\ \hline\hline
164 & 15 & 16 \\ \hline
165 & 15 & 16  \\ \hline
166 & 16 & 17 \\ \hline
167 & 17 & 18  \\ \hline
184 & 16 & 17  \\ \hline
188 & 17 & 18 \\ \hline
217 & 16 & 17  \\ \hline
220 & 16 & 17 \\ \hline
\end{tabular}
\end{center}
\bigskip

In these cases the bound $L_2(n)$ is equal or greater than $L(n)$, except for the
case $n=161$. 

The two functions $L(n)$ and $\Vert n\Vert$ coincide in 771 values  of $n$
in the range $1\le n\le 1000$,  the difference being equal to $1$ for the
229 other values in this range  with
a few exceptions.

\section{Particular values.}

\subsection{Numbers with small complexity.}
A good lower bound for $\Vert n\Vert$ is obtained from the knowledge of the 
largest number we may write with $m$ ones.  That is, given $m$, which is the 
largest natural number $N$ with $\Vert N\Vert=m$. The answer roughly is that
we must  group the $m$ ones in groups of three and multiply them.  To show this
we define the concept of \emph{extremal expression}. Let $\texttt{M}_m$ be an 
expression with $\Vert\texttt{M}_m\Vert=m$ (that is  $\texttt{M}_m$ is formed
with $m$ symbols $\texttt{x}$ and the operations of sum and product), and such
that its value $v(\texttt{M}_m)$ is the  maximum of all the expression with 
$m$ ones,
i.~e.
\begin{displaymath}
N=v(\texttt{M}_m)=\sup_{\Vert \texttt{A}\Vert =m} v(\texttt{A}).
\end{displaymath}
We say that such an  expression $\texttt{M}_m$ is extremal.

In the above situation $\Vert N\Vert =m$. In fact, since $N=v(\texttt{M}_m)$ and 
$\Vert \texttt{M}_m\Vert=m$, we have $\Vert N\Vert \le m$. Assume, by contradiction, 
that $\Vert N\Vert<m$. Then there will exists an expression \texttt{B} such that 
$v(\texttt{B})=N$ and $\Vert \texttt{B}\Vert=\Vert N\Vert<m$. Let $d$ be such that
$m=d+\Vert \texttt{B}\Vert$. We may construct an expression $\texttt{C}$ such that
$\texttt{C}=\texttt{B}+\texttt{x}+\cdots+\texttt{x}$ and such that $\Vert \texttt{C}
\Vert=\Vert \texttt{B}\Vert+d=m$ and $v(\texttt{C})=v(\texttt{B})+d>N$. This 
contradicts the definition of $\texttt{M}_m$.

It is easy to see that the following expressions are extremal
\begin{gather*}
\texttt{M}_1=\texttt{x},\quad \texttt{M}_2=\texttt{(x + x)},\quad 
\texttt{M}_3=\texttt{(x +
(x+x))},\\
\texttt{M}_4=\texttt{(x+x)(x+x)},
\quad \texttt{M}_5=\texttt{(x+(x+x))(x+x)},\dots
\end{gather*}
We see that given $m$ the extremal expression $\texttt{M}_m$ is not unique.
For example for $m=4$ the expression $\texttt{M}_4=\texttt{(x+(x+(x+x)))}$ 
is another possibility.

We shall use here  a  not very precise notation. For example, we shall write 
$\texttt{M}_3^a\texttt{M}_2$ to denote any expression having this form, not 
defining how the product is constructed from its factors. So, $\texttt{M}_3^4$
denotes any of the expressions $\texttt{((M}_3\texttt{M}_3\texttt{)(M}_3\texttt{M}_3
\texttt{))}$, $\texttt{(M}_3\texttt{(M}_3\texttt{(M}_3\texttt{M}_3\texttt{)))}$
or any other form of grouping the factors.

\begin{proposition}
Let $\texttt{M}_2=\texttt{(x + x)}$, $\texttt{M}_3=\texttt{(x +
(x+x))}$ and $\texttt{M}_4=\texttt{(x+x)(x+x)}$. For $n>1$, the expressions
\begin{displaymath}
\texttt{M}_n=\begin{cases}
\texttt{M}_3^{\,k} & \text{if $n=3k$},\\
\texttt{M}_3^{\,k-1}\texttt{M}_4 & \text{if $n=3k+1$},\\
\texttt{M}_3^{\,k}\texttt{M}_2 & \text{if $n=3k+2$},
\end{cases}
\end{displaymath}
are extremal. 
\end{proposition}

\begin{proof}
We may check the proposition for $n=2$, $3$ and $4$ directly.

Assume the assertion for all $s<n$ and try to prove it for $n\ge5$. Certainly there is
one extremal expression \texttt{K} with $\Vert\texttt{K}\Vert=n$. Then there are
two expressions \texttt{A} and \texttt{B} such that 
$\texttt{K}=\texttt{(A}+\texttt{B)}$ or $\texttt{K}=\texttt{(A}\texttt{B)}$. 
\texttt{A} and \texttt{B} are extremal expressions because \texttt{K} is extremal.
We may replace \texttt{A} and \texttt{B} by extremal expressions of the same complexity
and value and the resulting expression $\texttt{K}'$  will be also extremal.
Hence, without loss of generality, we may assume, using the induction hypothesis,
that \texttt{A} and \texttt{B}  are of the form given in the Proposition or 
$\texttt{A}=\texttt{x}$ and \texttt{B} is as in the Proposition.

The case $\texttt{K}=\texttt{(A}+\texttt{B)}$ it is only possible if $v(\texttt{A})$
or $v(\texttt{B})=1$, because, in other cases, the expression $\texttt{(AB)}$
contradicts the extremality of \texttt{K}.  But $\texttt{K}=\texttt{(x}+\texttt{M}_3^k
\texttt{)}$, $\texttt{K}=\texttt{(x}+\texttt{M}_3^{k-1}\texttt{M}_4\texttt{)}$,
or $\texttt{K}=\texttt{(x}+\texttt{M}_3^{k}\texttt{M}_2\texttt{)}$ are impossible with
$n\ge5$. Because these expressions are clearly not extremal. (Compare with
$\texttt{M}_3^{\,k-1}\texttt{M}_4$, $\texttt{M}_3^{\,k}\texttt{M}_2$ or 
$\texttt{M}_3^{\,k+1}$ respectively).

Therefore $\texttt{K}=\texttt{(A}\texttt{B)}$ where \texttt{A} and \texttt{B} 
are like those in the Proposition. Some of the combinations are not possible:
for example $\texttt{A}=\texttt{M}_3^{\,k}\texttt{M}_2$ and 
$\texttt{B}=\texttt{M}_3^{\,j-1}\texttt{M}_4$ are not possible since 
$\texttt{M}_3^{\,k+j-1}\texttt{M}_4\texttt{M}_2$ is improved by 
$\texttt{M}_3^{\,k+j+1}$ and \texttt{K} will not be extremal. A case by case analysis
proves that \texttt{K} is  one of the three forms in the Proposition.
\end{proof}

\begin{corollary}
For $a=0$, $1$, or $2$ and $b\in\N$ we have:
\begin{displaymath}
\Vert 2^a3^b\Vert=2a+3b, \qquad a=0, 1, 2.
\end{displaymath}
\end{corollary}

All natural numbers $n>1$ may be written in a unique way as $n=2a+3b$ with $a=0$, 
$1$ or $2$. In this case $2^a3^b$ is the greatest number $m$ with $\Vert m\Vert=n$.
Hence $m>2^a3^b$ implies $\Vert m\Vert>2a+3b$.

We define $g$ by
\begin{displaymath}
g(n)=\begin{cases} 3a & \text{if $n\in[3^a, 3^a+3^{a-1})$},\\
3a+1 & \text{if $n\in[3^a+3^{a-1},2\cdot3^a)$},\\
3a+2 & \text{if $n\in[2\cdot3^a,3^{a+1})$},
\end{cases}
\end{displaymath}
we  then have $g(n)\le \Vert n\Vert$ for each $n$. 
\begin{corollary}\label{antnueve}
For any $n\ge2$ we have
\begin{displaymath}
3\frac{\log n}{\log3}\le\Vert n\Vert \le L(n)\le3\frac{\log n}{\log2}.
\end{displaymath}
\end{corollary}

\begin{proof}
We only need to prove  the first inequality. If $n=3^a$, we see directly that
the inequality is true.  If $x\in(3^a, 3^a+3^{a-1}]$, we have $\Vert x\Vert\ge3a+1$.
Then
\begin{displaymath}
\Vert x\Vert\ge\Vert 3^a\Vert+1=3 a+1\ge
3\frac{\log (4\cdot 3^{a-1})}{ \log3}\ge3\frac{\log x}{\log3}.
\end{displaymath}
Analogously for $x\in(4\cdot3^{a-1},2\cdot 3^a]$ we have
\begin{displaymath}
\Vert x\Vert\ge \Vert 4\cdot 3^{a-1}\Vert+1\ge 3a+2\ge 3\frac{\log(2 \cdot
3^a)}{\log 3}.
\end{displaymath}
Finally for $x\in\bigl(2\cdot3^{a}, 3^{a+1}\bigr]$, we only need to check that
\begin{displaymath}
\Vert x\Vert\ge \Vert 2\cdot 3^a\Vert+1=3a+3\ge 3\frac{\log (3^{a+1})}{
\log 3}.
\end{displaymath}

\end{proof}

\section{The problem $\PNP$ and the complexity of the natural numbers.}

\subsection{Idea of the problem $\PNP$.}
Before explaining the problem we must describe  the classes $\mathbf{P}$
and $\mathbf{NP}$.  Consider a finite alphabet $A$, and let $A^*$ be the set of 
\emph{words}, that is, the set of finite sequences of elements of $A$.

We call \emph{language} a subset $S\subset A^*$. We say that $S$ is in the
class $\mathbf{P}$ if there is an algorithm $T$ and a polynomial $p(t)$ such that
with a word $x$ as input, $T$ gives an output $T(x)$, such that $T(x)=1$ if $x\in S$ 
and $T(x)=0$ if $x\notin S$. Also $T$ gives the output $T(x)$ in a time bounded by
$p(|x|)$ (here $|x|$ denotes the length of the word $x$). We then  say that $T$ is a 
polynomial algorithm. In a few words we may say that $\mathbf{P}$ is the class of
languages recognizable in polynomial time. It is important to notice that this
concept is very stable with respect to the diverse definitions of what is an 
algorithm, how we compute the ``time'' that the algorithm $T$ takes to give the
output, or even if we consider the same language in a different alphabet (as when
we consider a set of natural numbers written in different basis).  In other 
words,
the concept does not change if we give proper definitions of these concepts.

The class $\mathbf{NP}$ consists of the languages recognizables by non deterministic
polynomial algorithms. That is $S\subset A^*$ is in $\mathbf{NP}$ if there exists 
an algorithm T and a polynomial $p(x)$ such that for each $x\in S$ there is $y\in A^*$
with $|y|\le p(|x|)$ and such that with the input $(x,y)$ the algorithm gives 
the output $T(x,y)=1$ in time bounded by $p(|x|)$. On the other hand if $x\notin S$
we have $T(x,y)=0$ for all $y$ with $|y|\le p(|x|)$.

We say that in this case $T$ is a non-deterministic algorithm since  to 
obtain $x\in S$ we must first choose $y$. If we know which $y$ to take this 
process is fast, but if we do not know $y$, we may try each possible $y$, but this
will need a time $\ge |A|^{p(|x|)}$ which in practice is impossible.

Again the class $\mathbf{NP}$ is very stable with respect to possible changes in
the definitions. Also many practical problems are in this class.

It is easy to check that $\mathbf{P}\subset\mathbf{NP}$.  The question is whether
these two classes are the same. To understand a bit more of the difficulty
observe the following.

Our experience as mathematicians teaches us that to understand a proof, or better
to check the correctness of a proof is a task of type $\mathbf{P}$. That is the 
time needed is proportional to the length of the proof.

On the other hand to determine if a conjecture $x$ is a Theorem we need first to 
write the proof $y$ and then apply the above procedure to check the correctness of 
the pair $(x,y)$.  The set of Theorems is not in the class $\mathbf{NP}$ since as 
we know the length of the proof $|y|$  is not bounded by the length of the theorem $x$,
that is $|y|\not\le p(|x|)$. But for each polynomial $p(t)$, the following set is 
in $\mathbf{NP}$
\begin{displaymath}
\mathcal{T}_p=\{x: x \text{ is a theorem with a proof of length }\le p(|x|)\}.
\end{displaymath}
Maybe someone finds these definitions rather vague, but the formal logic allows
one to make things precise.

If $\mathbf{P}=\mathbf{NP}$ and the proof were sufficiently constructive (technically, that we can find a polynomial algorithm for an $\mathbf{NP}$-complete problem), 
then there
would exist a polynomial algorithm  that would allow not only decide if 
$x\in\mathcal{T}_p$, but also to find in this case a proof for $x$ in polynomial time.
The mathematicians would not be needed any more.

When one recalls the achievements of the $20^{\text{th}}$ century: proof of Fermat's theorem,
classification of finite simple groups, pointwise convergence of Fourier series 
of function in $L^p$, Riemann's hypothesis for algebraic varieties over fields of
characteristic $p$,  independence of continuum hypothesis,  and many more, one gets
the impression that there exists an algorithm to decide $x\in\mathcal{T}_p$, by
searching directly for a proof, not by trial and error. This algorithm consists in 
taking promising students, give them the possibility to travel and speak with
specialists on the topic in question, let them try to solve analogous questions, study the solution of related problems, and so on \dots

\section{Connection of the complexity of natural numbers and the problem $\PNP$.}
     
Consider the assertion $\Vert4787\Vert=28$. We may decompose it in two parts. The first,
$\Vert4787\Vert\le 28$, has a very easy proof
\begin{equation}\tag{$\star$}\label{star}
4787= 2 + 3(2+3^2)(1+2^43^2).
\end{equation}
The other part of the assertion $\Vert4787\Vert\ge 28$, has a much more laborious 
proof. Just now I do not know any other way than computing the values of 
$\Vert n\Vert$ for all $n\le4787$, a task that, on my personal computer, took several
hours.  

Of course this does not imply that it is easy to find proof as in \eqref{star}.

Consider the sets
\begin{displaymath}
A=\{(n,c)\in\N^2: \Vert n\Vert\le c\},\qquad
B=\{(n,c)\in\N^2: \Vert n\Vert> c\}.
\end{displaymath}
The fact, as we have remarked, that if $(n,c)\in A$, then there is a 
relatively short proof of it, shows us that $A$ is in the class $\mathbf{NP}$.

Roughly, a set $A$ is in $\mathbf{NP}$, if to prove that $x\in A$  
an exhaustive search is required, which in principle is exponential in the size of 
$x$, but once the proof has been found, it is easily recognized (in polynomial time 
with respect to the size of $x$).
Complete information may be found in the book \cite{GJ}. These problems 
bring to mind the one of finding a needle in a haystack. Once we have found 
the needle there
is no doubt that the task is done, but at first it appears unreachable  since the 
straw
is so similar to the needle that we do not see 
any other means than search methodically.

The core of the problem $\PNP$ is whether in situations where there exists a
short proof, there is always 
a direct path to find it.  If $\mathbf{P}=\mathbf{NP}$,
then there is always a direct path  to the proof without hesitations.  At first 
sight this appears a wild assumption, but the rigorous proof of 
$\mathbf{P}\ne\mathbf{NP}$ eludes us still after twenty seven years of study. 

Recently Microsoft has funded an investigation center and has contracted 
Michael Friedman, (Fields medal  in 1986). Friedman has the intention of trying 
to solve the question $\PNP$. Microsoft  will invest 2.6 million 
dollars each year in this program. 

It appears that $\mathbf{P}=\mathbf{NP}$ is false, but not all is so simple. Sometimes
tasks that appear to need an exhaustive search have been proved simple. We shall give
an example.

Let $\mathcal{C}\subset\N$ the set of composite numbers. At first sight it appears
that the only means to proof that $n$ is composite is to divide $n$ by each number
$m\le \sqrt{n}$ and check if some remainder equals $0$. The size of $n$ is of the order
of the number of digits needed to write it, i.e. of the order $\log n$. The number
of needed checks maybe $\sqrt{n}=e^{(\log n)/2}$, which grows exponentially with
$\log n$. And if really $n$ is composite there is a short proof: to exhibit a 
proper divisor $d$ of $n$. That is $\mathcal{C}$ is in the class $\textbf{NP}$.

But it is not so difficult to decide whether $n$ is composite. If $n$ is prime and 
$b$ is prime with $n$ we have $b^{n-1}\equiv1\bmod{n}$. An idea somewhat more 
elaborate, let $n$ be a prime and $n-1=2^s t$,  in the sequence of the rests of
$b^t$, $b^{2t}, \dots, b^{2^s t}$ $\bmod n$  the last different from $1$ must be
$-1$. In the other case it is sure that $n$ is composite.
This is the famous Miller-Rabin test. It is known that if the generalized Riemann
hypothesis is true, then if $n$ is composite, the test of Miller-Rabin  is not 
satisfied for some $b<2(\log n)^2$. Hence, under the mentioned hypothesis, we 
have a fast algorithm (polynomial)  to decide whether $n$ is composite: to do the 
test of Miller-Rabin for all $b<2(\log n)^2$. 

Another incentive to pose the problem $\PNP$ is the existence of 
\emph{$\mathbf{NP}$-complete problems}. That is sets $B\subset \N$  such that 
$B$ is in the class $\mathbf{NP}$ and, for which from  $B\in\mathbf{P}$ it follows 
that $\mathbf{P}=\mathbf{NP}$. 

From Euclid's times,  mathematicians have had a clear concept of 
 algorithm. Turing gives a further step and by an effort of introspection gives us
a precise definition.  Turing's mental image is that of a mathematician, notebook in hand, computing.  By abstracting the procedure Turing created the idea of a modern 
computer. Starting from Turing's definitions it is possible to quantify  the time
a computer will spend on a given task and so to give a precise definition of the classes
$\mathbf{P}$ and $\mathbf{NP}$.

The first connection of the complexity of the natural numbers with the problem 
$\PNP$ is the fact that $\mathbf{P}=\mathbf{NP}$ implies the existence of a fast
algorithm to compute $\Vert n\Vert$. There will be constants $C$ and $k\in \N$ and 
an algorithm  that will compute $\Vert n\Vert$ in  time $\le C(\log n)^k$. 

\section{Complexity of boolean functions.}
There is another connection, this time structural, between the complexity of natural 
numbers and the problem $\PNP$. To explain this connection we must define a 
related concept, that of the complexity of a boolean function.

The set $\{0,1\}$ is a field when we consider the composition laws sum and product
$\bmod\ 2$. For each number $n$ let $\mathcal{F}_n$ be the set of functions $f\colon
\{0,1\}^n\to\{0,1\}$. The set $\mathcal{F}_n$ is a ring if we take sum and product 
with respect to the field in the image $\{0,1\}$.

\def\1{{\bold1}}\def\0{{\bold 0}}
For example consider the constant functions $\1$, $\0$ and the components
${\boldsymbol\pi}_j$  defined by
${\boldsymbol \pi}_j({\bold x})={\boldsymbol \pi}_j(x_1, x_2, \dots, x_n)=x_j$. 

The ring $\mathcal{F}_n$ is generated by these functions, i.~e.~we may write
any function $f\in\mathcal{F}_n$ as a polynomial of the above functions. To see this
given ${\boldsymbol\varepsilon}=(\varepsilon_1, \dots,
\varepsilon_n)\in\{0,1\}^n$, we define the function $f_{\boldsymbol
\varepsilon}=\prod_j(\delta_j+{\boldsymbol
\pi}_j)$, where, for each $j$, $\delta_j=1+\varepsilon_j$. Then 
$f_{\boldsymbol\varepsilon}({\bold x})=0$, except for ${\bold x}={\boldsymbol\varepsilon}$. Hence, any function $g$ may be written
\begin{displaymath}
g=\sum_{{\boldsymbol\varepsilon}\in S}f_{\boldsymbol\varepsilon},
\end{displaymath}
where $S$ is the set of ${\boldsymbol\varepsilon}$ such that
$g({\boldsymbol\varepsilon})=1$. 

As in the case of the natural numbers, we may define the complexity of the 
elements of $\mathcal{F}_n$.  It will be the greatest function  $f\mapsto\Vert f\Vert$
such that
\begin{displaymath}
\Vert\0\Vert=\Vert\1\Vert=0;\quad \Vert {\boldsymbol\pi}_j\Vert=1;\qquad 
\Vert f+g\Vert\le \Vert f\Vert+\Vert g\Vert;\quad
\Vert fg\Vert\le \Vert f\Vert+\Vert g\Vert.
\end{displaymath}

For any $\theta\in(0,1)$, most of the elements of ${\mathcal F}_n$ have complexity
$\ge2^{\theta n}$. The proof of this result is done by counting how many 
elements have complexity $k$, say  $a_k$. It is easy to see that $a_0=2$, $a_1=2n$.
From $f$ and $g$ with $\Vert f\Vert=j$ and $\Vert g\Vert=k-j$ we get, at most, four
elements with complexity $\le k$. They are $f+g$, $fg$, $1+f+g$, $1+fg$. 
With these observations we get 
\begin{displaymath}
a_k\le 4(a_1a_{k-1}+a_2a_{k-2}+\cdots+a_{k-1}a_1).
\end{displaymath}
It follows that $a_k\le A_k$, where $A_k$ is defined by
\begin{displaymath}
A_0=2;\quad A_1=2n;\quad A_k=4\sum_{j=1}^{k-1}A_jA_{k-j}.
\end{displaymath}
From this definition we get 
\begin{displaymath}
\sum_{k=0}^\infty A_kx^k=\frac{17-\sqrt{1-32nx}}{8}; \qquad
A_k=\frac{1}{2(2k-2)}{2k-2\choose k}(8n)^k.
\end{displaymath}
Hence
\begin{displaymath}
a_k\le A_k\sim\frac{2^{5k}}{8\sqrt{2\pi}k^{3/2}}n^k.
\end{displaymath}
Therefore for $x$ large
\begin{displaymath}
\sum_{k=0}^x A_k\le c\sum_{k=0}^x(32 n)^k\le c'(32n)^x\le Ae^{Bx\log n},
\end{displaymath}
hence if $x<2^{\theta n}$, with $0<\theta<1$, we get
\begin{displaymath}
\sum_{k=0}^x A_k\ll\card(\mathcal{F}_n)=2^{2^n},
\end{displaymath}
proving our assertion.

Each construction of $f(x_1, \dots, x_n)$ as a polynomial allows one to prove
an assertion of type $\Vert f\Vert\le a$. But from the polynomial expression
we may get something more practical: a circuit that allows to compute 
$f(x_1, \dots, x_n)$  starting from the imputs $x_j$. 

As in the case of natural numbers, it is difficult to prove inequalities of 
type  $\Vert f\Vert> a$. In fact the situation is surprising: we have seen 
that in the set of functions with $n$ variables, the complexity is 
usually larger than $2^{\theta n}$. Hence one would expect to have an 
easy task in defining 
a sequence of functions $(f_n)$, where $f_n$ depends on $n$
variables and such that  $\Vert f_n\Vert>2^{\theta n}$. On the contrary 
it has only been achieved that $\Vert f_n\Vert>p(n)$, where $p$ is a polynomial of 
small degree (see \cite{Z}, \cite{H}).  The problem here is not to prove
that there exist sequences with $\Vert f_n\Vert>2^{\theta n}$, which, as we have
seen is easy, but to define explicitly a  concrete sequence of functions 
for which this is so. When we speak of  ``define explicitly'' we refer to 
a technical concept that needs some explanation. We must exclude easy solutions
as: \emph{let $f_n$ the first function of $n$ variables with maximum complexity}.
We say that $(f_n)$ is given explicitly if there is an algorithm that computes
the value of $f_n(x_1,\dots,x_n)$ in a reasonable time.

The problem $\PNP$ induces one to consider a special sequence of boolean functions.
Let $a$ be a natural number and consider $n=\binom{a}{2}$ the number of pairs.
Our variables will be
\begin{displaymath}
x_{12},\, x_{13},\, x_{23},\, x_{14},\, x_{24},\, x_{34}, 
\,\dots,\,x_{1a},\, x_{2a},\, \cdots,\,
x_{a-1\, a}.
\end{displaymath}
In this way, each set of values of these variables $\in\{0,1\}^n$ may be seen 
as a graph with $a$ vertices and where $x_{jk}=1$ if and only if the vertices 
$j$ and $k$ are connected by an edge of the graph. For each $b\le a$ let 
$f_b^a(x_{12},\dots, x_{a-1\,a})$ be the function that is equal $1$ if and only if
there is a set of $b$ vertices  such that all of then are connected in the graph.

It is plausible that $\Vert f_b^a\Vert\ge\binom{a}{b}$, since  to compute the 
value of $f_b^a$ in a given graph we need to check each set of $b$ vertices.
It can be shown that, if this is so, then $\mathbf{P}\ne \mathbf{NP}$. In this way
to prove $\Vert f_b^a\Vert\ge\binom{a}{b}$ is, I think, the most promising path 
to solve the $\PNP$ question.
\bigskip

In the case of the complexity of natural numbers, an analogous question 
is the following, posed by Guy \cite{G}:

\begin{problem} Is there a sequence of natural numbers $(a_n)$ such that
\begin{equation}\tag{1}\label{uno}
\lim_{n\to\infty}\frac{\Vert a_n\Vert}{\log a_n}>\frac{3}{\log 3}?
\end{equation}
\end{problem}
A good candidate is the sequence $2^n$. All computed values  satisfy 
$\Vert 2^n\Vert=2n$. Selfridge asks (see \cite{G}) whether there exists 
any $n$ with 
$\Vert 2^n\Vert<2n$.

If for some $n$ and $k$ we would have $2^n=3^k$, (which is clearly impossible), 
the second expression would give us  $\Vert 2^n\Vert<2n$. Of course the advantage
would be greater for big $n$ than for small $n$.  Although the above is impossible,
maybe another type of equality would yield $\Vert 2^n\Vert<2n$. For example, if 
for some $n$,  $2^n$ written in base $3$ has small digits. Again, this is unlikely
but not impossible. Also, there may exist another type of expression 
of $2^n$. The question here is whether a number  of the form
\begin{displaymath}
(1+1)(1+1)\cdots(1+1),
\end{displaymath}
may be written in some way with fewer $1$'s.  We have almost a trivial example
$4=(1+1)(1+1)=1+1+1+1$. 
Here we have the same number of $1$'s so that I call it an almost-example. Maybe 
there are non-trivial almost-examples, for example
\begin{displaymath}
2^{27} = 1 +(1+2\cdot 3)(1+2^3\cdot 3^2)(1+2^9\cdot 3^3(1+2\cdot 3^2)).
\end{displaymath}
If we replace each $2$ by $1+1$ and each $3$ by $1+1+1$ we get an expression 
for $2^{27}$ with $57$ ones, in which the multiplicative structure of $2^{27}$ 
is not used.

The above equality proves that $\Vert2^{27}-1\Vert\le 56$. In spite of an intense 
search I have not found an $n>2$ such that $\Vert 2^n-1\Vert<2n-1$, but I think 
this may happen. 

The evidence appears to be in favor of the existence of a sequence that satisfies 
\eqref{uno}. For example, we may look at figure 2. There we have put a little
disk with center at each point $(n,\Vert n\Vert)$ with $1\le n\le 2000$ and also 
we have drawn the smooth curves that bound $\Vert n\Vert$, i.~e.~$3(\log t)/\log 3$ 
and $3(\log t)/\log 2$, and also the curve $5\log t/2\log2$. The points overlap and 
we see some lines parallel to the $x$-axis. We see that the upper bound appears to 
be bad and that apparently $\Vert n\Vert\le 5\log n/2\log2$, whereas in reality we 
have only proved that this inequality is true for almost all $n\in\N$.

    \begin{figure}[H]
    \begin{center}
    \includegraphics[width = 0.6\textwidth]{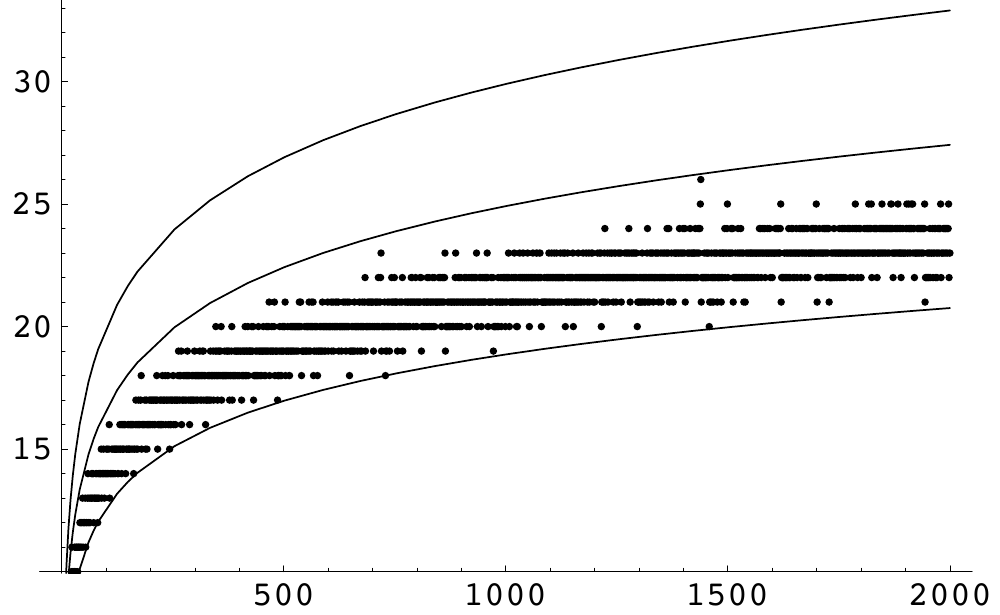}
   \caption{Graph of $\Vert n\Vert$.}
    \end{center}
    \label{tres}
    \end{figure} 
    
But this figure says nothing about the limit $\lim\Vert n\Vert/\log n$, 
in which we are interested in. We only see that for the first 2000 values
of $n$ this sequence is bounded by  the limits $5/2\log 2$ and $3/\log 3$. 

\section{Conjectures}

I have computed, using Proposition \ref{computing}, the complexity of the 
first $200\,000$ natural numbers. Looking at these numbers, one sees many regularities.
We will call them conjectures about the behavior of the function $\Vert\cdot\Vert$, 
although I have not much confidence in that they persist for larger numbers.

These conjectures were derived from tables such as this one
$$\begin{matrix} 3 & 6 & 9 & 12 & 15 &18 & 21 & 24 \\
\noalign{\medskip\hrule\medskip}
 10   & 100   & 1000   & 10000   & 100000
      & 1000000   & 10000000   & 100000000
      \\    & 22   & 220   & 2200   & 22000
      & 220000   & 2200000   & 22000000   \\ 
      & 21   & 210   & 2101   & 21010   & 210100
      & 2101000   & 21010000   \\    & 
      & 202   & 2100   & 21000   & 210000   & 
   2100000   & 21000000   \\    &    & 201
      & 2020   & 20200   & 202000   & 2020000
      & 20200000   \\    &    & \overline{122}   & 2010
      & 20100   & 201000   & 2010000   & 20100000
      \\    &    &    & 2002   & 20020
      & 200222   & 2002220   & 20022200   \\ 
      &    &    & 2001   & 20010   & 200200
      & 2002000   & 20020000   \\    & 
      &   & \overline{1221 }  & 20002   & 200100   & 
   2001000   & 20010000   \\    &    & 
      & 1220   & 20001   & 200020   & 2000200
      & 20002000   \\    &   & 
      & 1212   & \overline{12221 }  & 200010   & 2000100
      & 20001000   \\    &    & 
      & 1211   & 12210   & 200002   & 2000020
      & 20000200   \\    &    & 
      & 1201   & 12200   & 200001   & 2000010
      & 20000100   \\    &    & 
      & 1122   & 12122   & \overline{122210 }  & 2000002
      & 20000020   \\    &   & 
      & 1121   & 12120   & 122100   & 2000001
      & 20000010   \\    &    & 
      & 1112   & 12111   & 122000   & \overline{1222100}
      & 20000002   \\    &    & 
      &    & 12110   & 121220   & 1221000
      & 20000001   \\   &   & 
      &    & 12102   & 121200   & 1220000
      & \overline{12221000 }  \\   &    & 
      &    & 12101   & 121121   & 1212200
      & 12210000   \\    &    & 
      &    & 12012   & 121110   & 1212000
      & 12200000   \\   &    & 
      &   & 12010   & 121100   & 1211210
      & 12122000   \\   &    & 
      &    & 12001   & 121022   & 1211100
      & 12121201   \\   &    & 
      &   & 11221   & 121020   & 1211000
      & 12120000   \\ \end{matrix} $$

In this table we have written in columns the numbers with complexity $3n$ ($n=1$, $2$, \dots, $8$), written in base $3$ and in decreasing order.

The first observation: $\Vert 3n\Vert=3+\Vert n\Vert$ is wrong. $\Vert107\Vert=16$ and 
$\Vert 321\Vert=\Vert1+2^6 5\Vert=18$. But  the following
conjectures seem to be true:

\begin{conjecture} For each natural number $n$, there is an integer $a\ge0$ such that
$\Vert 3^j n\Vert = 3(j-a)+\Vert 3^an\Vert$ for each natural number $j\ge a$.
\end{conjecture}

Let us define the set $A=\{n\in\N : 
\Vert 3^j n\Vert =3j+\Vert n\Vert \text{ for all } j\}$. 

\begin{conjecture}
For each pair of natural numbers $p$ and $q$, there exists $a\ge0$ such 
that, for $j\ge a$, we have $\Vert p(q3^j +1)\Vert = 3j+1+\Vert p\Vert+
\Vert q\Vert$. 
\end{conjecture}

The main observation in the above table is that the greatest numbers with
complexity $3n$ are those natural numbers contained in the sequence 
$(3^n a_n)$, where $a_n$ is given by
\begin{multline*}
1, {2(3+1)\over 3^2}, {2^6\over 3^4}, {2\cdot 3
+1\over3^2}, {2( 3^2+1)\over 3^3}, {2\cdot 3^2+1\over 3^3},
{2^9\over 3^6},\\ {2( 3^3+1)\over 3^4}, {2\cdot 3^3+1\over 3^4},
\dots, {2( 3^k+1)\over 3^{k+1}}, {2\cdot 3^k+1\over 3^{k+1}},\dots 
\end{multline*}

\begin{conjecture}
There exist three transfinite sequences $(a_\alpha)_{\alpha<\xi}$, 
$(b_\alpha)_{\alpha<\xi}$,
$(c_\alpha)_{\alpha<\xi}$ of rational numbers, such that the (greatest)
numbers of complexity $3n$ (respectively $3n+1$, $3n+2$) are the (first) 
natural numbers contained in the sequence $(3^n a_\alpha)$, 
(resp.~$(3^n b_\alpha)$, $(3^n c_\alpha)$).

$\xi$ is an infinite numerable ordinal such that $\omega\xi=\xi$. 
\end{conjecture}

These sequences start in the following way:
\def\hu{\phantom{3}}
\begingroup
\Small
\begin{align*}
&(a_\alpha), &  
&1, \hu{8\over 9}, \hu{64\over81},\hu {7\over9}, {20\over 27},\,\dots \to{2\over3} &  
&{160\over243}, {52\over81},\,\dots \to{16\over27}
 & {1280\over2187}, {140\over243},\dots \to{5\over9}\,\dots\\
&(b_\alpha), &  &{4\over3},
{32\over27}, {10\over9}, {256\over243}, {28\over27},\,\dots \to{1} &
&{80\over81},{26\over27},\,\dots \hu\to{8\over9} & {640\over729},
{70\over81},\dots\to{64\over81}\,\dots\\  
&(c_\alpha),& 
&2,{16\over9}, {5\over3}, {128\over81}, {14\over9},\,\dots\to{4\over3} & 
&{320\over243}, {35\over27},\,\dots\to{32\over27} 
&{95\over 81}, {2560\over2187},\dots\to{10\over9}\,\dots 
\end{align*}
\endgroup
where the dots indicate infinite sequences, and where the indicated limits
are not terms of  the sequences.

\begin{conjecture}
The three sequences are decreasing. The denominators of each term $a_\alpha$,
$b_\alpha$ or $c_\alpha$ are powers of $3$.
\end{conjecture}

\begin{conjecture}
The numbers of the sequence $(a_\alpha)$ are the numbers of the set
\begin{displaymath}
\left\{\frac{n}{3^{\Vert n\Vert/3}}: \Vert n\Vert\equiv0 \mod3, \quad
\text{and}\quad n\in A\right\},
\end{displaymath}
ordered decreasingly. 
\end{conjecture}

\begin{conjecture}
The numbers of the sequence $(b_\alpha)$ are the numbers of the set
\begin{displaymath}
\left\{\frac{n}{3^{(\Vert n\Vert-1)/3}}: \Vert n\Vert\equiv1 \mod3,\quad
\text{and}\quad n\in A\right\},\end{displaymath}
ordered decreasingly. 
\end{conjecture}

\begin{conjecture}
The numbers of the sequence $(c_\alpha)$ are the numbers of the set
\begin{displaymath}
\left\{\frac{n}{3^{(\Vert n\Vert-2)/3}}: \Vert n\Vert\equiv2 \mod3, \quad
\text{and}\quad n\in A\right\},\end{displaymath}
ordered decreasingly. 
\end{conjecture}

The following conjectures are more doubtful. They are only based on a few
cases.

\begin{conjecture}
For all ordinals $\beta<\xi$ we have
\begin{displaymath}
\lim_{n\to\infty}a_{\omega\beta+n}=c_\beta/3, \quad
\lim_{n\to\infty}b_{\omega\beta+n}=a_\beta,\quad
\lim_{n\to\infty}c_{\omega\beta+n}=b_\beta.
\end{displaymath}
\end{conjecture}
This is the basis of the assertion about the value of $\xi$, which appears to be
at least $\xi=\omega^\omega$, since this is the least solution of $\omega\xi=\xi$.

The following assertions, along with conjecture 8, allow to predict, with some accuracy, the values
of the transfinite sequences.

\begin{conjecture}
The numbers of the sequence $b_{\omega\beta+n}$ that converges to $a_\beta=b/3^a$
(with $\Vert b\Vert=3a$) are numbers from the sequences
\begin{displaymath}
\frac{p(q 3^j+1)}{3^{a+j}}, \quad \hbox{\rm where } \quad
b=pq, \hbox{ and,  }\Vert p(q3^j+1)\Vert=3a+3j+1,  
\end{displaymath}
and those sporadic terms of the sequence $2^{3j+2}/3^{2j+1}$ contained between
$\sup_{\gamma<\beta}a_\gamma$ and $a_\beta$. 
\end{conjecture}

\begin{conjecture}
The numbers of the sequence $c_{\omega\beta+n}$ that converges to $b_\beta=b/3^a$
(with $\Vert b\Vert=3a+1$) are numbers from the sequences
\begin{displaymath}
\frac{p(q 3^j+1)}{3^{a+j}}, \quad \hbox{\rm where } \quad
b=pq, \hbox{ and,  }\Vert p(q3^j+1)\Vert=3a+3j+2,  
\end{displaymath}
and those sporadic terms of the sequence $2^{3j+1}/3^{2j}$ contained between
$\sup_{\gamma<\beta}b_\gamma$ and $b_\beta$. 
\end{conjecture}

\begin{conjecture}
The numbers of the sequence $a_{\omega\beta+n}$ that converges to $c_\beta/3=b/3^a$
(with $\Vert b\Vert=3a-1$) are numbers from the sequences
\begin{displaymath}
\frac{p(q 3^j+1)}{3^{a+j}}, \quad \hbox{\rm where } \quad
b=pq, \hbox{ and,  }\Vert p(q3^j+1)\Vert=3a+3j,  
\end{displaymath}
and those sporadic terms of the sequence $2^{3j}/3^{2j}$ contained between
$\sup_{\gamma<\beta}\frac13 c_\gamma$ and $\frac13 c_\beta$. 
\end{conjecture}

In Conjecture 9, 10 and 11 we observe that some terms come from subsequent 
sequences. For example, the term $c_\omega=320/243$ is the term corresponding
to $j=0$ of the sequence $2^6(4\cdot3^j+1)/3^{j+5}$, that converges to $b_3=256/243$. 

The above conjectures allow one to predict, for example, the 200 largest numbers 
with complexity $30$.

The numbers with complexity 14 divided by $81$, are 
\begin{align*}
c_0&={162\over81},&
c_1&={144\over81},&
c_2&={135\over81},&
c_3&={128\over81},\\
c_4&={126\over81},&
c_5&={120\over81},&
c_6&={117\over81},&
c_7&={114\over81},\\
c_9&={112\over81},&
c_{10}&={111\over81},&
c_{11}&={110\over81},&
c_{13}&={109\over81},\\
c_{\omega+1}&={105\over81},&
c_{\omega+2}&={104\over81},&
c_{\omega+3}&={102\over81},&
c_{\omega+6}&={100\over81},\\
c_{\omega+8}&={99\over81},&
c_{\omega+10}&={98\over81},&
c_{\omega+14}&={97\over81},&
c_{\omega2}&={95\over81},\\
c_{\omega2+3}&={93\over81},&
c_{\omega2+5}&={92\over81},&
c_{\omega2+8}&={91\over81},&
c_{\omega3+4}&={88\over81},\\
c_{\omega3+7}&={87\over81},&
c_{\omega3+15}&={86\over81},&
c_{\omega4+2}&={85\over81},&
c_{\omega5+1}&={83\over81},\\
c_{\omega^2+\omega+2}&={79\over81},&
c_{\omega^2+\omega2+3}&={77\over81},&
\end{align*}
\begin{displaymath}
{71\over81},\quad {69\over81},\quad {67\over81},\quad {59\over81},
\end{displaymath}
For the last four numbers I do not have enough data to know the corresponding 
ordinal. 




\end{document}